\documentclass[11pt, reqno]{amsart}
\usepackage[dvipsnames,usenames]{color}
\usepackage{hyperref}
\usepackage{graphicx}
\usepackage{epsfig}
\usepackage[latin1]{inputenc}
\usepackage{amsmath}
\usepackage{amsfonts}
\usepackage{amssymb}
\usepackage{amsthm}
\usepackage{amscd}
\usepackage{verbatim}
\usepackage{subfigure}
\usepackage{caption}
\usepackage{pinlabel}
\usepackage{stmaryrd}
\usepackage{enumerate, enumitem}
\usepackage{todonotes}
\usepackage{bm}
\usepackage{thmtools}
\usepackage{thm-restate}
\usepackage{lipsum}
\usepackage{setspace}
\usepackage{mathtools}
\usepackage[all]{xypic}
\usepackage[abs]{overpic}
\usepackage{color}
\usepackage[normalem]{ulem}

\allowdisplaybreaks

\usepackage{tikz}
\usetikzlibrary{arrows}
\usetikzlibrary{decorations.pathreplacing}
\usepackage{verbatim}
\usetikzlibrary{cd}
\tikzset{taar/.style={double, double equal sign distance, -implies}}
\tikzset{amar/.style={->, dotted}}
\tikzset{dmar/.style={->, dashed}}
\tikzset{aar/.style={->, very thick}}

\newtheorem{theorem}{Theorem}[section]

\newtheorem{conjecture}[theorem]{Conjecture}
\newtheorem{question}[theorem]{Question}

\theoremstyle{definition}
\newtheorem{definition}[theorem]{Definition}

\theoremstyle{remark}
\newtheorem{remark}[theorem]{Remark}
\newtheorem{example}[theorem]{Example}

\newcommand{\alphas}{\boldsymbol{\alpha}}
\newcommand{\betas}{\boldsymbol{\beta}}

\def\F{\mathbb{F}}
\def\N{\mathbb{N}}
\def\Q{\mathbb{Q}}
\def\R{\mathbb{R}}
\def\Z{\mathbb{Z}}

\def\lebang{<^!}

\def\cC{\mathcal{C}}

\def\cH{\mathcal{H}}

\def \gr {\operatorname{gr}}

\def\d{\partial}

\def\co{\colon}

\def\spinc{\textrm{Spin}^c}

\def\id{\textup{id}}

\def\CF {\mathit{CF}}
\def\HF {\mathit{HF}}
\newcommand\HFhat{\widehat{\HF}}

\newcommand\HFp {\HF^+}

\newcommand \CFm {\CF^-}
\newcommand \HFm {\HF^-}

\def\CFK{\mathit{CFK}}
\def\HFK{\mathit{HFK}}

\newcommand\HFKm{\mathit{HFK}^-}
\newcommand\HFKhat{\widehat{\mathit{HFK}}}

\def\spinc {{\operatorname{spin^c}}}

\newcommand{\skewsimeq}{\mathrel{\rotatebox[origin=c]{-180}{$\simeq$}}}

\newcommand{\CFgroup}{\mathfrak{D}}
\newcommand{\CFKgroup}{\mathfrak{C}}
\newcommand{\CFKmgroup}{\mathfrak{C}'}
\newcommand{\CFIgroup}{\mathfrak{I}}
\newcommand{\CFIgrouphat}{\hat{\mathfrak{I}}}
\newcommand{\SF}{\textit{SF}}
\newcommand{\CFKIgrouphat}{\hat{\mathfrak{I}}_K}

\author[J.\ Hom]{Jennifer Hom}
\thanks{The author was partially supported by NSF grants DMS-1552285 and DMS-2104144.}
\address {School of Mathematics, Georgia Institute of Technology, Atlanta, GA 30332}
\email{hom@math.gatech.edu}

\numberwithin{equation}{section}

\title[Homology cobordism, knot concordance, and HF homology]{Homology cobordism, knot concordance, and Heegaard Floer homology}

\begin{document}

\begin{abstract} 
We review some recent results in knot concordance and homology cobordism. The proofs rely on various forms of Heegaard Floer homology. We also discuss related open problems.
\end{abstract}

\maketitle

\section{Introduction}\label{sec:intro}
Many interesting phenomena occur in 3- and 4-dimensions that do not occur in higher dimensions. Indeed, the Poincar\'e conjecture in dimensions five and higher was proved by Smale \cite{Smale} in 1961 using techniques from surgery theory, while the Poincar\'e conjecture in dimension three remained unsolved for another 40 years until Perelman's proof of Thurston's geometrization conjecture \cite{Perelman1, Perelman2, Perelman3}. In 1982, Freedman \cite{Freedman} proved the topological 4-dimensional Poincar\'e conjecture, while the smooth 4-dimensional Poincar\'e conjecture remains open.
As another example of how dimension four is special, by work of Freedman \cite{Freedman} and Donaldson  \cite{DonaldsonBAMS}, $\R^n$ admits smooth structures that are not diffeomorphic to the standard one only when $n = 4$.\footnote{Smale, Perelman, Freedman, and Donaldson all won Fields medals for their work discussed here; Perelman declined the award.}

The 3-dimensional homology cobordism group and the knot concordance group are fundamental structures in low-dimensional topology. The former played a key role in Manolescu's disproof of the high dimensional triangulation conjecture \cite{ManolescuPin2Triangulation}, while the latter has the potential to shed light on the smooth 4-dimensional Poincar\'e conjecture (see, for example, \cite{FGMW, ManolescuPiccirilllo}).


Our goal is to review some recent applications of Heegaard Floer theory to homology cobordism and knot concordance, and to discuss the power and limitations of these tools to address major open questions in the field.

\subsection{Homology cobordism} Two closed, oriented 3-manifolds $Y_0, Y_1$ are \emph{homology cobordant} if there exists a smooth, compact, oriented 4-manifold $W$ such that $\d W = -Y_0 \sqcup Y_1$ and the inclusions $\iota \co Y_i \to W$ induce isomorphisms
\[ \iota_* \co H_*(Y_i; \Z) \to H_*(W; \Z) \]
for $i=0,1$. The key point is that, on the level of homology, $W$ looks like a product. The \emph{3-dimensional homology cobordism group} $\Theta^3_\Z$ consists of integer homology 3-spheres modulo homology cobordism, under the operation induced by connected sum. A homology sphere $Y$ represents the identity in $\Theta^3_\Z$ if and only if $Y$ bounds a homology 4-ball, and the inverse of $[Y]$ in $\Theta^3_\Z$ is $[-Y]$, where $-Y$ denotes $Y$ with the opposite orientation. The Rokhlin invariant \cite{Rokhlin} gives a surjective homomorphism
\[ \mu \co \Theta^3_\Z \to \Z/\Z2, \]
showing that $\Theta^3_\Z$ is nontrivial. Manolescu \cite{ManolescuPin2Triangulation} showed that if $\mu(Y)=1$, then $Y$ is not of order two in $\Theta^3_\Z$. By work of Galewski-Stern \cite{GalewskiStern} and Matumoto \cite{Matumoto}, this leads to a disproof of the triangulation conjecture in dimensions $\geq 5$. See \cite{ManolescuLectTri, ManolescuICM} for an overview of this work. The triangulation conjecture is also false in dimension 4, by work of Casson; see \cite{AkbulutMcCarthy}.

Fintushel-Stern \cite{FintushelStern85} used gauge theory to show that $\Theta^3_\Z$ is infinite, and Furuta \cite{furuta} and Fintushel-Stern \cite{FintushelStern90} improved this result to show that $\Theta^3_\Z$ contains a subgroup isomorphic to $\Z^\infty$. Fr\o yshov \cite{Froyshov} used Yang-Mills theory to define a surjective homomorphism $\Theta^3_\Z \to \Z$, showing that $\Theta^3_\Z$ has a direct summand isomorphic to $\Z$. (This is stronger than having a $\Z$ subgroup, since, for example, $\Z$ is a subgroup of $\Q$ but not a summand.) In joint work with Dai, Stoffregen, and Truong, we use Hendricks-Manolescu's involutive Heegaard Floer homology \cite{HMInvolutive} to prove the following:

\begin{theorem}[{\cite{DHSThomcob}}]\label{thm:DHSThomcob}
The homology cobordism group $\Theta^3_\Z$ contains a direct summand isomorphic to $\Z^\infty$.
\end{theorem}

Fundamental questions about the structure of $\Theta^3_\Z$ remain open:

\begin{question}\label{quest:thetatorsion}
Does $\Theta^3_\Z$ contain any torsion? Modulo torsion, is $\Theta^3_\Z$ free abelian?
\end{question}

If there is any torsion in $\Theta^3_\Z$, two-torsion seems the most likely. Indeed, any integer homology sphere $Y$ admitting an orientation-reversing self-diffeomorphism is of order at most two in $\Theta^3_\Z$. However, as far as we are aware, all known examples of such $Y$ bound integer homology balls, and hence are trivial in $\Theta^3_\Z$.

In a different direction, it is natural to ask which types of manifolds can represent a given class $[Y] \in \Theta^3_\Z$. The first answers to this question were in the positive. Livingston \cite{livingston} showed that every class in $\Theta^3_\Z$ can be represented by an irreducible integer homology sphere and Myers \cite{myers} improved this to show that every class has a hyperbolic representative.

In the negative direction, Fr\o yshov (in unpublished work), F. Lin \cite{LinExact}, and Stoffregen \cite{stoffregen-sums} showed that there are classes in $\Theta^3_\Z$ that do not admit Seifert fibered representatives. Nozaki-Sato-Taniguchi \cite{NST} improved this result to show that there are classes that do not admit a Seifert fibered representative nor a representative that is surgery on a knot in $S^3$. However, none of these results were sufficient to obstruct $\Theta^3_\Z$ from being generated by Seifert fibered spaces. In joint work with Hendricks, Stoffregen, and Zemke, we prove the following:

\begin{theorem}[{\cite{HHSZ, HHSZ2}}]\label{thm:SFS}
The homology cobordism group $\Theta^3_\Z$ is not generated by Seifert fibered spaces. More specifically, let $\Theta_{\SF}$ denote the subgroup generated by Seifert fibered spaces. The quotient  $\Theta^3_\Z/\Theta_{\SF}$ is infinitely generated.
\end{theorem}

In light of the aforementioned Nozaki-Sato-Taniguchi result, it is natural to ask:

\begin{question}\label{quest:surgeries}
Do surgeries on knots in $S^3$ generate $\Theta^3_\Z$?
\end{question}

\noindent The expectation is that surgeries on knots in $S^3$ are not sufficiently generic to generate $\Theta^3_\Z$, but such a result seems beyond the capabilities of current tools.

\subsection{Knot concordance} Two knots $K_0, K_1 \subset S^3$ are \emph{concordant} if there exists a smooth, properly embedded annulus $A$ in $S^3 \times [0,1]$ such that $K_i = A \cap (S^3 \times \{i\})$ for $i=0,1$. The \emph{knot concordance group} $\cC$ consists of knots in $S^3$ modulo concordance, under the operation induced by connected sum. The inverse of $[K]$ in $\cC$ is given by $[-K]$, where $-K$ denotes the reverse of the mirror image of $K$. A knot $K \subset S^3 = \partial B^4$ is \emph{slice} if it bounds a smoothly embedded disk in $B^4$. Fox-Milnor \cite{FoxMilnor} and Murasugi \cite{Murasugi} showed that $\cC$ is nontrivial, and J. Levine \cite{JLevine} used the Seifert form to define a surjective homomorphism $\cC \to \Z^\infty \oplus (\Z/2\Z)^\infty \oplus (\Z/4\Z)^\infty$, demonstrating that $\cC$ is in fact highly nontrivial. In higher odd dimensions (that is, knotted $S^{2n+1}$ in $S^{2n+3}$, $n \geq 1$), Levine's homomorphism is an isomorphism, while in the classical dimension, the kernel is nontrivial \cite{CassonGordon}. 
See \cite{Livingstonsurvey} for a survey of knot concordance.

We can consider various generalizations of the knot concordance group. For example, rather than considering annuli in $S^3 \times [0,1]$, we may consider annuli in homology cobordisms. 
Two knots $K_0 \subset Y_0$ and $K_1 \subset Y_1$ are \emph{homology concordant} if they cobound a smooth, properly embedded annulus in a homology cobordism between $Y_0$ and $Y_1$.

Let $\cC_\Z$ denote the group of knots in $S^3$, modulo homology concordance. A knot $K \subset S^3$ represents the identity in $\cC_\Z$ if and only if $K$ bounds a smoothly embedded disk in some homology 4-ball. Since $B^4$ is of course a homology 4-ball, there is naturally a surjection from $\cC$ to $\cC_\Z$. A natural question is whether or not this map is injective; in other words:

\begin{question}\label{quest:ZHBslice}
If a knot $K \subset S^3$ bounds a disk in a homology 4-ball, must $K$ also bound a disk in $B^4$?
\end{question}

One reason why the question above is challenging is that many obstructions to a knot $K$ bounding a disk in $B^4$ also obstruct $K$ from bounding a disk in a homology 4-ball.

An even more difficult question is the following:

\begin{question}\label{quest:S4PC}
If a knot $K \subset S^3$ bounds a disk in a homotopy 4-ball, must $K$ also bound a disk in $B^4$?
\end{question}

Recall the smooth 4-dimensional Poincar\'e conjecture which by work of Freedman \cite{Freedman} may be stated as follows:

\begin{conjecture}[Smooth 4-dimensional Poincar\'e conjecture]\label{conj:S4PC}
If a smooth 4-manifold $X$ is homeomorphic to $S^4$, then $X$ is actually diffeomorphic to $S^4$.
\end{conjecture}

A negative answer to Question \ref{quest:S4PC} provides one possible strategy for disproving Conjecture \ref{conj:S4PC}.  Indeed, by Freedman \cite{Freedman}, any homotopy 4-sphere is homeomorphic to $S^4$. Now suppose we found a homotopy 4-sphere $X$ and a knot $K \subset S^3 = \d (X \setminus \mathring{B}^4)$ such that $K$ bounds a smoothly embedded disk in $X \setminus \mathring{B}^4$. If we could obstruct $K$ from bounding a smoothly embedded disk in $B^4$, then it follows that $X$ cannot be diffeomorphic to $S^4$. This approach was attempted in \cite{FGMW, ManolescuPiccirilllo}, but has yet to lead to a disproof of the conjecture.


We now return to the group $\cC_\Z$. This group is naturally a subgroup of $\widehat{\cC}_\Z$, the group of manifold-knot pairs $(Y, K)$, where $Y$ is a homology sphere bounding a homology ball and $K$ is a knot in $Y$, modulo homology concordance. One can ask whether the injection from $\cC_\Z$ to $\widehat{\cC}_\Z$ is a surjection. Adam Levine \cite{Levinenonsurj} answered this question in the negative, showing that there exist knots in a homology null-bordant $Y$ (in fact, his example bounds a contractible 4-manifold) that are not concordant to any knot in $S^3$. Expanding on this result, in joint work with Levine and Lidman, we prove the following:

\begin{theorem}[{\cite{HLL}}]\label{thm:HLL}
The subgroup $\cC_\Z \subset \widehat{\cC}_\Z$ is of infinite index. More specifically, 
\begin{enumerate}
	\item the quotient $\widehat{\cC}_\Z/\cC_\Z$ is infinitely generated,  and 
	\item \label{it:HLL2} contains a subgroup isomorphic to $\Z$.
\end{enumerate}
\end{theorem}

\noindent This result demonstrates the vast difference between knots in $S^3$ and knots in arbitrary homology spheres bounding homology balls, up to concordance. The examples in the infinite generation part of the theorem bound contractible 4-manifolds; it's unknown whether the examples in the $\Z$ subgroup do. Zhou \cite{Zhou} proved that the quotient $\widehat{\cC}_\Z/\cC_\Z$ has a subgroup isomorphic to $\Z^\infty$. (It's unknown whether his examples bound contractible 4-manifolds.) Forthcoming joint work with Dai, Stoffregen, and Truong \cite{DHSThomconc} improves Zhou's result to a $\Z^\infty$-summand.

One can also consider concordance in more general 4-manifolds. Let $R$ be a ring. Two closed, oriented, connected 3-manifolds $Y_0, Y_1$ are \emph{$R$-homology cobordant} if there exists a smooth, compact, oriented 4-manifold $W$ such that $\d W = -Y_0 \sqcup Y_1$ and the inclusions $\iota \co Y_i \to W$ induce isomorphisms
\[ \iota_* \co H_*(Y_i; R) \to H_*(W; R) \]
for $i=0,1$. We have already discussed the case $R=\Z$. The rational homology cobordism group $\Theta^3_\Q$ contains elements of order two, for example, $[\R P^3]$; in contrast, as asked in Question \ref{quest:thetatorsion}, it remains open whether there is any torsion in the integer homology cobordism group $\Theta^3_\Z$.

We can consider concordances in other $R$-homology cobordisms, such as $\Q$-homology cobordisms. 
A knot $K \subset S^3$ is \emph{rationally slice} if it is $\Q$-homology concordant to the unknot, or equivalently, if $K$ bounds a smoothly embedded disk in a rational homology 4-ball. 

Let $\cC_{\Q S}$ denote the subgroup of $\cC$ consisting of rationally slice knots. Cochran, based on work of Fintushel-Stern \cite{Fintushel-Stern:1984-1}, showed that the figure-eight knot is rationally slice. Hence $\Z/2\Z$ is a subgroup of $\cC_{\Q S}$, since the figure-eight is negative amphichiral and not slice. Cha \cite{Cha:2007-1} extended this result to show that $\cC_{\Q S}$ has a subgroup isomorphic to $(\Z/2\Z)^\infty$. A natural question to ask is whether $\cC_{\Q S}$ contains elements of infinite order (see, for example, \cite[Problem 1.11]{GeorgiaConference}). Joint work with Kang, Park, and Stoffregen uses the involutive knot Floer package of Hendricks-Manolescu \cite{HMInvolutive} to prove:

\begin{theorem}[{\cite{HKPS}}]\label{thm:Qslice}
The group of rationally slice knots $\cC_{\Q S}$ contains a subgroup isomorphic to $\Z^\infty$.
\end{theorem}

The figure-eight is slice in a rational homology 4-ball $W$ with $H_1(W; \Z) = \Z/2\Z$ (see, for example, \cite[Section 3]{AkbulutLarson}), as are Cha's examples \cite{Cha:2007-1}.

\begin{question}\label{quest:QHBslice}
Does there exist a knot $K \subset S^3$ that is not slice in $B^4$ but is slice in a rational homology 4-ball $W$ with $|H_1(W; \Z)|$ odd?
\end{question}

\noindent Compare this question to Question \ref{quest:ZHBslice}, which asks whether there is a knot $K \subset S^3$ that is not slice in $B^4$ but is slice in a integer homology 4-ball $W$. Indeed, both Questions \ref{quest:ZHBslice} and \ref{quest:QHBslice} can be viewed as incremental steps towards Question \ref{quest:S4PC}, a negative answer to which would in turn disprove the smooth 4-dimensional Poincar\'e conjecture.

\subsection{Ribbon concordance}
\begin{figure}[htb!]
\includegraphics[scale=0.8]{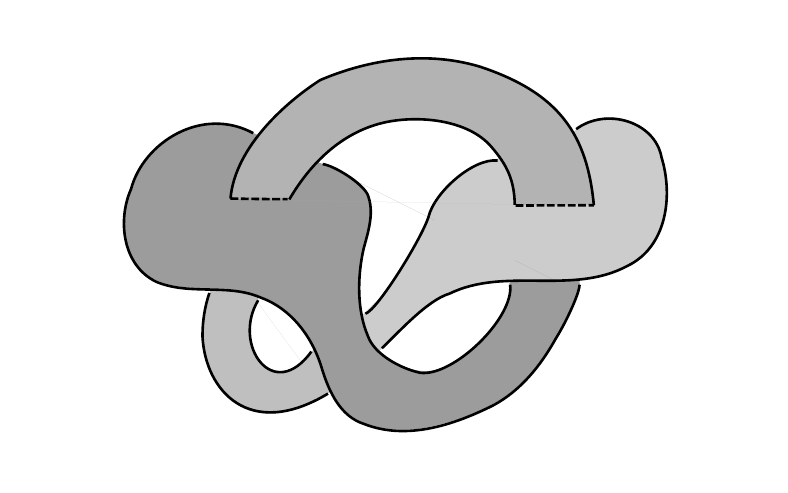}
\caption{\label{fig:ribbon} An example of a ribbon disk.}
\end{figure}

We conclude the introduction with a discussion of ribbon knots, ribbon concordances, and ribbon homology cobordisms. A knot $K \subset S^3$ is \emph{ribbon} if it bounds an immersed disk in $S^3$ with only ribbon singularities. A \emph{ribbon singularity} is a closed arc consisting of intersection points of the disk with itself such that the preimage of this arc is two disjoint arcs in the disk with one arc $a_1$ being contained entirely in the interior of the disk and the other arc $a_2$ having its endpoints on the boundary of the disk. See Figure \ref{fig:ribbon}. Note that ribbon knots are slice, since ribbon singularities can be resolved in the 4-ball (namely, by pushing the arc $a_1$ farther into the 4-ball).

\begin{conjecture}[Slice-ribbon conjecture \cite{FoxProblems}]
Every slice knot is ribbon.
\end{conjecture}

\noindent The slice-ribbon conjecture is true for two-bridge knots \cite{Lisca} and many infinite families of pretzel knots \cite{GreeneJabuka, Lecuona}. On the other hand, potential counterexamples exist; see, for example, \cite{GST, AbeTagami}.

Equivalently, a ribbon knot can be defined as a knot $K \subset S^3$ that bounds a \emph{ribbon disk} in $B^4$, that is, a smoothly embedded disk $D \subset B^4$ such that the radial Morse function on $B^4$ restricted to $D$ has no interior local maxima. There is a \emph{ribbon concordance} from $K_0$ to $K_1$ if there is a concordance from $K_0$ to $K_1$ in $S^3 \times [0,1]$ with no interior local maxima (with respect to the natural height function on $S^3 \times [0,1]$), or equivalently if projection to $[0,1]$ is Morse with only index 0 and index 1 critical points. Note that, unlike ordinary concordance, ribbon concordance is not symmetric (and that the convention regarding the direction of the ribbon concordance varies in the literature).

\begin{conjecture}[\cite{Gordon}]\label{conj:GordonRibbon}
Ribbon concordance is a partial order. That is, if there exists a ribbon concordance from $K_0$ to $K_1$ and a ribbon concordance from $K_1$ to $K_0$, then $K_0 = K_1$.
\end{conjecture}

\noindent Gordon \cite{Gordon} proved that Conjecture \ref{conj:GordonRibbon} holds for fibered knots and two-bridge knots, and more generally for the class of knots generated by such knots under the operations of connected sum and cabling. He also proved that if $S$ is a ribbon concordance from $K_0$ to $K_1$, then $\pi_1(S^3 \setminus \nu (K_0)) \to \pi_1(X)$ is injective and $\pi_1(S^3 \setminus \nu (K_1)) \to \pi_1(X)$ is surjective, where $X$ denotes the exterior of $S$ in $S^3 \times [0,1]$.

The notion of a ribbon concordance can also be generalized to homology cobordisms. A \emph{ribbon cobordism} between two 3-manifolds is a cobordism admitting a handle decomposition with only 1- and 2-handles. Observe that the complement of a ribbon concordance from one knot to another is naturally a ribbon cobordism between their knot complements. We have the following 3-manifold analog of Conjecture \ref{conj:GordonRibbon}:

\begin{conjecture}[\cite{DLVW}]\label{conj:ribbon3}
Ribbon $\Q$-homology cobordism is a partial order on closed, oriented, connected 3-manifolds. That is, if there exists a ribbon $\Q$-homology cobordism from $Y_0$ to $Y_1$ and a ribbon $\Q$-homology cobordism from $Y_1$ to $Y_0$, then $Y_0$ and $Y_1$ are diffeomorphic.
\end{conjecture}

\noindent The results in \cite{DLVW} provide evidence in support of Conjecture \ref{conj:ribbon3}.

\subsection{Organization}
The remainder of this article is devoted to discussing applications of Heegaard Floer homology to the theorems and problems discussed above. We will describe various Heegaard Floer chain complexes associated to 3-manifolds and knots inside of them. As we progress, our chain complexes will have more and more structure; we will sketch how this additional structure leads to the results described in the introduction.

In Section \ref{sec:HF}, we discuss properties of the Heegaard Floer 3-manifold invariant of \cite{OSDisks, OSProperties} and applications to homology cobordism. In Section \ref{sec:HFK}, we move on to knot Floer homology  \cite{OSKnots, RasmussenKnots} and applications to concordance. With the advent of involutive Heegaard Floer homology  \cite{HMInvolutive}, these invariants can be endowed with additional structure; in Section \ref{sec:HFI}, we describe involutive Heegaard Floer homology, and in Section \ref{sec:involutiveCFK}, we delve into involutive knot Floer homology. Lastly, in Section \ref{sec:future}, we discuss the potential (or lack thereof) for Heegaard Floer homology to answer the questions posed in the Section \ref{sec:intro}.

For an introduction to many of the tools described in this article, we refer the reader to Sections 1--3 of \cite{HomPCMI}.

\section*{Acknowledgements}
I would like to thank my collaborators Irving Dai, Kristen Hendricks, Sungkyung Kang, Adam Levine, Tye Lidman, JungHwan Park, Matthew Stoffregen, Linh Truong, and Ian Zemke for the privilege of working with them and for their patience in working with me. The questions posed in this article are generally well-known in the field; I learned of most of them through problems sessions and discussions over the years. I am grateful to all of my colleagues who have generously given their time and energy to conference organization and conversation. I would also like to thank Robert Lipshitz, Chuck Livingston, Ciprian Manolescu, and Andr\'as Stipsicz for helpful comments on an earlier draft.

\section{Heegaard Floer homology: the 3-manifold invariant}\label{sec:HF}

In this section, we consider the Heegaard Floer 3-manifold invariant and maps induced by 4-dimensional cobordisms. For expository overviews of Heegaard Floer homology, see, for example, \cite{OSintro}, \cite{JSurvey}, and \cite[Section 2]{HomPCMI}.

\subsection{Properties and examples}
Given a closed, oriented 3-manifold $Y$, its Heegaard Floer homology $\HFm(Y)$ is a finitely generated, graded module over $\F[U]$, where $\F = \Z/2\Z$ and $U$ is a formal variable in degree $-2$. (There are other flavors, $\HFp(Y), \HFhat(Y)$ of Heegaard Floer homology, but for the purposes of this article, we will focus on the minus version.) 

More precisely, every closed, oriented 3-manifold $Y$ can be described as a union of two handlebodies; such a decomposition is called a \emph{Heegaard splitting}. In turn, a Heegaard splitting can be described via a \emph{Heegaard diagram} $\cH$, consisting of a closed, oriented surface $\Sigma$ of genus $g$, together with $g$ $\alpha$-circles and $g$ $\beta$-circles, which describe how the handlebodies fill in the surface on either side. (These circles are required to satisfy a certain homological condition.) For technical reasons, we also fix a basepoint $z$ in the complement of the $\alpha$- and $\beta$-circles. Any two diagrams representing the same 3-manifold can be related by a sequence of \emph{Heegaard moves}, as described in \cite[Section 2.6]{OSintro}; see also \cite[Section 1]{HomPCMI}.

From this data, Ozsv\'ath-Szab\'o \cite{OSDisks} construct $\CFm(\cH)$, a free, finitely generated, graded chain complex over $\F[U]$. The variable $U$ keeps track of the basepoint $z$. The chain homotopy type of $\CFm(\cH)$ is an invariant of $Y$; that is, it does not depend on the choice of Heegaard diagram, nor on any other choices made in the construction. We often write $\CFm(Y)$ to denote this chain homotopy class, or a representative thereof. Juh\'asz-Thurston-Zemke \cite{JTNaturality} prove something even stronger: Heegaard Floer homology is natural, in the sense that it assigns a concrete module, rather than an isomorphism class of modules, to a 3-manifold.

\begin{example}
The Heegaard Floer homology of $S^3$ is $\HFm(S^3) = \F_{(0)}[U]$, where the subscript $(0)$ denotes that $1 \in \F[U]$ is in grading $0$. (This can easily be computed using the definition of the Heegaard Floer chain complex.)
\end{example}

\begin{example}
The Heegaard Floer homology of the Brieskorn homology sphere $\Sigma(2,3,7)$ is $\HFm(\Sigma(2,3,7)) = \F_{(0)}[U] \oplus \F_{(0)}$. (This is not so easy to compute directly from the definition of Heegaard Floer homology; however, it is a straightforward consequence of some of the formal properties of Heegaard Floer homology.)
\end{example}

\begin{example}
The Heegaard Floer homology of the Brieskorn homology sphere $\Sigma(2,3,5)$ is $\HFm(\Sigma(2,3,5)) = \F_{(-2)}[U]$. (This can be computed using some of the formal properties of Heegaard Floer homology.)
\end{example}

\begin{remark}
In this article, we take the slightly unconventional grading convention above, which simplifies the formula for gradings in, for example, $\HFm$ of connected sums. Many other sources use the convention that $\HFm(S^3) = \F_{(-2)}[U]$, which simplifies calculations in $\HFp$.
\end{remark}

\begin{remark}
For rational homology spheres, the gradings in Heegaard Floer homology take values in $\Q$. For integer homology spheres, the gradings take values in $2\Z$. For 3-manifolds $Y$ with $H_1(Y; \Z)$ infinite, the gradings are slightly more complicated; see \cite[Section 4.2]{OSIntersectionForms}.
\end{remark}

Since the degree of $U$ is $-2$, any homogenously graded polynomial in $\F[U]$ is of the form $U^n$ for some $n \in \N$. Thus, by the fundamental theorem of finitely generated graded modules over a PID, we have that $\HFm(Y)$ is of the form
\[ \bigoplus_{i=1}^N \F_{(d_i)}[U] \oplus \bigoplus_{j=1}^M \F_{(c_j)}[U]/U^{n_j}\F[U], \]
for $Y$ a rational homology sphere; that is, $\HFm(Y)$ is a direct sum of a free part and a $U$-torsion part.
Ozsv\'ath-Szab\'o \cite[Theorem 10.1]{OSProperties} show that when $Y$ is an integer homology sphere, $N=1$, that is, $\HFm(Y)$ is of the form
\[ \HFm(Y) = \F_{(d)}[U] \oplus \bigoplus_{j=1}^M \F_{(c_j)}[U]/U^{n_j}\F[U]. \]
The number $d$ above is called the \emph{$d$-invariant} of $Y$, denoted $d(Y)$. More generally, when $Y$ is a rational homology sphere with $|H_1(Y; \Z)|=k$, there are exactly $k$ free summands in $\HFm(Y)$, in which case one obtains a $k$-tuple of $d$-invariants of $Y$.

Heegaard Floer homology satisfies a K\"unneth-type formula under connected sums \cite[Theorem 1.5]{OSProperties}. That is, $\CFm(Y_1 \# Y_2)$ is chain homotopy equivalent to
\[ \CFm(Y_1) \otimes_{\F[U]} \CFm(Y_2). \]
In particular, if $Y_1$ and $Y_2$ are integer homology spheres, then the $d$-invariant is additive under connected sum. (An analogous statement also holds for more general 3-manifolds.)

\subsection{Cobordism maps} Heegaard Floer homology is a $(3+1)$-dimensional topological quantum field theory (TQFT) \cite[Theorem 1.1]{OSTriangles}. That is, to a 3-manifold, Heegaard Floer homology associates a module, and to a 4-manifold cobordism $W$ from $Y_0$ to $Y_1$, it associates a chain map
\[ F_W \co \CFm(Y_0) \to \CFm(Y_1). \]
When $W$ is a homology cobordism, $F_W$ induces an isomorphism
\[ (F_W \otimes \id)_* \co U^{-1} \HFm(Y_0) \to U^{-1} \HFm(Y_1), \]
where $U^{-1} \HFm(Y) = H_*(\CFm(Y) \otimes_{\F[U]} \F[U, U^{-1}])$ \cite[Proof of Theorem 9.1]{OSIntersectionForms}. A straightforward algebra calculation then implies that $d$-invariants are invariants of homology cobordism. In particular, we have a homomorphism
\[ d \co \Theta^3_\Z \to 2\Z, \]
and this homomorphism is surjective, since $d(\Sigma(2,3,5)) = -2$.

In light of the discussion above, and motivated by the desire to study the homology cobordism group $\Theta^3_\Z$, one could define an equivalence relation $\sim$, called \emph{local equivalence}, on Heegaard Floer chain complexes, where $\CFm(Y_0) \sim \CFm(Y_1)$ if there exist $\F[U]$-module chain maps
\[ f \co \CFm(Y_0) \to \CFm(Y_1) \qquad \textup{ and } \qquad g \co \CFm(Y_1) \to \CFm(Y_0) \]
inducing isomorphisms on $U^{-1}\HFm(Y_i)$. We can now consider the group
\[  \CFgroup = \{ \CFm(Y) \mid Y \textup{ an integer homology sphere} \} / \sim \]
under the operation induced by tensor product. This construction yields a homomorphism 
\[ \Theta^3_\Z \to  \CFgroup \]
obtained by sending $[Y]$ to $[\CFm(Y)]$. However, it turns out that $\CFgroup$ is isomorphic to $\Z$, with the isomorphism being given by $[\CFm(Y)] \mapsto d(Y)/2$.

\subsection{Ribbon homology cobordisms}
Ribbon homology cobordisms induce particularly nice maps on Heegaard Floer homology:
\begin{theorem}[{\cite[Theorem 1.19]{DLVW}}]\label{thm:DLVW}
Let $W$ be a ribbon homology cobordism from $Y_0$ to $Y_1$. Then
\[ F_W \co \HFm(Y_0) \to \HFm(Y_1) \]
is injective, and includes $\HFm(Y_0)$ into $\HFm(Y_1)$ as a direct summand.
\end{theorem}

\noindent The proof of Theorem \ref{thm:DLVW} relies on considering the \emph{double $D(W)$} of $W$, formed by gluing $W$ to $-W$ along $Y_1$; they also prove that the analogous statement holds for ribbon $\Z/2\Z$-homology cobordisms. Their approach was inspired by \cite{ZemkeRibbon}.

\section{Knot Floer homology}\label{sec:HFK}
In this section, we will discuss the Heegaard Floer knot invariant, maps induced by concordances, and various concordance invariants arising from the knot Floer complex.
For expository overviews of knot Floer homology, see \cite[Section 10]{OSintro}, \cite{ManolescuKnotIntro}, \cite{Homsurvey}, \cite[Section 2]{HeddenWatson}, and \cite[Section 3]{HomPCMI}. Note that \cite{HomPCMI} uses the same notation and conventions used here (i.e., viewing the knot Floer complex as a module over a two-variable polynomial ring), while the others use a different but equivalent formulation in terms of filtered chain complexes.

\subsection{Properties and examples}
For simplicity, we will focus on knots in integer homology spheres. Let $K$ be a knot in an integer homology sphere $Y$. We can describe the pair $(Y, K)$ via a \emph{doubly pointed Heegaard diagram $\mathcal{H}$}, which consists of a Heegaard diagram for $Y$ with an extra basepoint $w$. The knot $K$ is the union of two arcs, specified by connecting the basepoint $w$ to the basepoint $z$ in the complement of the $\alpha$-arcs, pushed slightly into one handlebody, and connecting $z$ to $w$ in the complement of the $\beta$-arcs, pushed slightly into the other handlebody. See, for example, \cite[Section 1]{HomPCMI} for more details.

From this data, Ozsv\'ath-Szabo \cite{OSDisks} and independently J. Rasmussen \cite{RasmussenKnots} construct a chain complex $\CFK(\cH)$. One way of constructing this chain complex (see, for example, \cite[Section 1.5]{ZemAbsoluteGradings}) is as a free, finitely generated, bigraded chain complex over $\F[U, V]$, the second formal variable $V$ corresponding to the second basepoint $w$. As one would hope, the chain homotopy type of $\CFK(\cH)$ is an invariant of the pair $(Y, K)$, and does not depend on the choice of diagram, nor on any of the other choices made in the construction. We often write $\CFK(Y, K)$, or simply $\CFK(K)$ when $Y=S^3$, to denote this chain homotopy class, or a representative thereof. Moreover, like the 3-manifold version, knot Floer homology is natural \cite{JTNaturality}.

\begin{example}
The knot Floer complex of the unknot in $S^3$ is generated over $\F[U, V]$ by a single generator $x$ in bigrading $(0,0)$ with trivial differential. (This can be computed directly from the definition of the knot Floer chain complex.)
\end{example}

\begin{example}
The knot Floer complex of the right handed trefoil is generated over $\F[U, V]$ by $a, b,$ and $c$ with the following differentials and bigradings:
\begin{center}
\begin{tabular}{*{16}{@{\hspace{10pt}}c}}
\hline
&  && $\d$ &&  $\gr$ \\
\hline
& $a$ && $0$ && $(0, -2$) \\ 
& $b$ && $U a + V c$ && $(-1, -1)$ \\ 
& $c$ && $0$ && $(-2, 0)$ \\ 
\hline
\end{tabular}
\end{center}
(This can be computed directly from the definition of the knot Floer chain complex.)
\end{example}

There is not a simple characterization of finitely-generated, graded modules over $\F[U, V]$, since $\F[U, V]$ is not a PID. One way to obtain a module over a PID is to set $V=0$ on the chain level. (There is a symmetry between $U$ and $V$, so one could instead choose to set $U=0$.) Taking homology of the resulting chain complex, we obtain a version of knot Floer homology, namely $\HFKm(Y, K)$:
\[  \HFKm(Y,K) = H_*(\CFK(Y,K)/(V=0)). \]
If we prefer an even simpler algebraic structure, we can set both $U$ and $V$ equal to zero on the chain level. Taking homology of the resulting chain complex, we obtain another version of knot Floer homology, denoted $\HFKhat(Y, K)$:
\[ \HFKhat(Y, K) = H_*(\CFK(Y,K)/(U=V=0)). \]
Using a suitable renormalized bigrading, this is the version of knot Floer homology whose graded Euler characteristic is the Alexander polynomial.

Like the 3-manifold invariant, the knot Floer complex satisfies a K\"unneth-type formula under connected sums \cite[Theorem 1.5]{OSProperties}. That is, $\CFK(Y_1 \# Y_2, K_1 \# K_2)$ is chain homotopy equivalent to 
\[ \CFK(Y_1, K_1) \otimes_{\F[U,V]} \CFK(Y_2, K_2). \]

\subsection{Maps induced by concordances}
Knot Floer homology also behaves nicely under cobordisms. Consider a cobordism $(W, S)$ from $(Y_0, K_0)$ to $(Y_1, K_1)$; that is, $W$ is a 4-manifold cobordism from $Y_1$ to $Y_2$ and $S \subset W$ a properly embedded connected surface with boundary $-K_0 \sqcup K_1$. The pair $(W, S)$ induces a module homomorphism
\[ F_{W,S} \co \CFK(Y_0, K_0) \to \CFK(Y_1, K_1). \]
When $W$ is a homology cobordism and $S$ is an annnulus, $F_{W,S}$ induces an isomorphism 
\[ (F_{W,S} \otimes \id)_* \co (U, V)^{-1}\HFK(Y_0, K_0) \xrightarrow[]{\cong} (U, V)^{-1}\HFK(Y_1, K_1), \]
where $(U, V)^{-1}\HFK(Y_0, K_0) = H_*(\CFK(Y,K) \otimes_{\F[U,V]} \F[U, U^{-1}, V, V^{-1}]$ \cite[Theorem 1.7]{ZemAbsoluteGradings}. When $W= S^3 \times [0,1]$, we may simply write $F_S$ instead of $F_{W,S}$.

Using this additional structure, we define an equivalence relation on these chain complexes that is well-suited to studying the knot concordance group, and more generally, concordances in homology cobordisms. 

\begin{definition}
Two knot Floer complexes, $\CFK(Y_0, K_0)$ and $ \CFK(Y_1, K_1)$ are \emph{locally equivalent}, denoted $\CFK(Y_0, K_0) \sim \CFK(Y_1, K_1)$, if there exist $\F[U, V]$-module chain maps
\[ f \co \CFK(Y_0, K_0) \to \CFK(Y_1, K_1) \qquad \textup{ and } \qquad g \co \CFK(Y_1, K_1) \to \CFK(Y_0, K_0) \]
inducing isomorphisms on $(U, V)^{-1} \HFK(Y_i, K_i)$. 
\end{definition}

\begin{remark}
In the literature, this equivalence relation is also referred to as $\nu^+$-equivalence \cite{KimPark} and stable equivalence \cite{Homsurvey}.
\end{remark}

As in the 3-manifold case above, we can now consider the group
\[ \CFKgroup = \{ \CFK(Y, K) \mid Y \textup{ a $\Z H S^3$ bounding a $\Z H B^4$, } K \textup{ a knot in $Y$} \} / \sim \]
under the operation induced by tensor product, giving us a homomorphism
\[  \widehat{\cC}_\Z \to \CFKgroup \]
obtained by sending $[(Y, K)]$ to $[\CFK(Y, K)]$. (We require $Y$ to bound a $\Z H B^4$ to parallel the definition of $\widehat{\cC}_\Z$.) By precomposing with the map $\cC \to \widehat{\cC}_\Z$, we obtain a homomorphism $\cC \to \CFKgroup$.

The group $\CFKgroup$ is not easy to study. One way to obtain a simpler algebraic structure is to set $V=0$, in which case we can run the analogous construction; that is, we can consider the group
\[ \CFKmgroup = \{ \CFK(Y, K)/(V=0) \mid Y \textup{ a $\Z H S^3$ bounding a $\Z H B^4$, } K \textup{ a knot in $Y$} \} / \sim \]
where $\CFK(Y_0, K_0)/(V=0) \sim \CFK(Y_1, K_1)/(V=0)$ if there exist $\F[U]$-module chain maps 
\begin{align*}
	f \co \CFK(Y_0, K_0)/(V=0) \to \CFK(Y_1, K_1)/(V=0) \\
	g \co \CFK(Y_1, K_1)/(V=0) \to \CFK(Y_0, K_0)/(V=0)
\end{align*}
inducing isomorphisms on $U^{-1} \HFKm(Y_i, K_i)$; we call this equivalence relation \emph{local equivalence mod $V$}. As in the 3-manifold case, this group $\CFKmgroup$ is isomorphic to $\Z$; indeed, up to renormalization, this construction yields the Ozsv\'ath-Szab\'o $\tau$-invariant \cite{OS4ballgenus} (see also \cite[Appendix A]{OST}).

In a case of mathematical Goldilocks, working over the full ring $\F[U,V]$ yields a group that is too complicated to study, 
while working over the ring $\F[U] = \F[U,V]/(V=0)$ yields a group that is too simple. Somewhat miraculously, it turns out that working over the ring $\F[U, V]/(UV=0)$ is just right, at least for knots in $S^3$. The main idea is that although $\F[U, V]/(UV=0)$ has zero-divisors, it is somehow closer to being a PID than $\F[U, V]$ is. Furthermore, for knots in $S^3$, the local equivalence group mod $UV$ is totally ordered, as we now describe.

\begin{definition}
We say that 
\[ \CFK(K_1) \leq \CFK(K_2) \]
if there exists an $\F[U, V]/(UV=0)$-module chain map 
\[ f \co \CFK(K_1) \to \CFK(K_2) \]
such that $f$ induces an isomorphism on $H_*(\CFK(K_i)/U)/V\textup{-torsion}$. If 
\[ \CFK(K_1) \leq \CFK(K_2) \quad \textup{ and } \quad \CFK(K_2) \leq \CFK(K_1), \]
then we say that $\CFK(K_1)$ and $\CFK(K_2)$ are \emph{locally equivalent mod $UV$}.
\end{definition}

\begin{remark}
Note that since $U$ and $V$ are both zero-divisors in $\F[U, V]/(UV)$, we cannot invert them directly. Also note that requiring $f$ to induce an isomorphism on $H_*(\CFK(K_i)/U)/V\textup{-torsion}$ is ever so slightly stronger than requiring $f$ to induce an isomorphism on $V^{-1} H_*(\CFK(K_i)/U)$.
\end{remark}

This is the same total order as the one induced by $\varepsilon$ in \cite{HomEpsilon}. Indeed,  one way to define the $\{-1, 0, 1\}$-valued concordance invariant $\varepsilon$ \cite{Homcables} (see also \cite[Section 3]{DHSTmore}) is as follows:
\begin{itemize}
	\item $\varepsilon(K) \leq 0 $ if and only if $\CFK(K) \leq \F[U, V]$,
	\item $\varepsilon(K) \geq 0$ if and only if $\CFK(K) \geq \F[U, V]$.
\end{itemize}
In particular, $\varepsilon(K) = 0$ if and only if $\F[U, V] \leq \CFK(K) \leq \F[U, V]$.

For knots in $S^3$, we are able to characterize their knot Floer complexes, up to local equivalence mod $UV$:

\begin{theorem}[{\cite[Theorem 1.3]{DHSTmore}}]\label{thm:phi}
Let $K$ be a knot in $S^3$. The knot Floer complex of $K$ is locally equivalent mod $UV$ to a standard complex, which can be represented by a finite sequence of nonzero integers. Moreover, if we endow the integers with the following unusual order
\[ -1 \lebang -2 \lebang -3 \lebang \dots \lebang 0 \lebang \dots \lebang 3 \lebang 2 \lebang 1 \]
then local equivalence classes mod $UV$ are ordered lexicographically with respect to their standard representatives.
\end{theorem}

\noindent See Section 4 of \cite{DHSTmore} for the definition of a standard complex. This characterization of knot Floer complexes up to local equivalence mod $UV$ is a key step in the definition of the linearly independent family of concordance homomorphisms 
\[ \varphi_i \co \cC \to \Z, \quad i \in \N \]
from \cite{DHSTmore}. The main idea is that $\varphi_i(K)$ is the signed count of the number of times that $i$ appears in the sequence of integers parametrizing the local equivalence class mod $UV$ of $\CFK(K)$. (For symmetry reasons, we actually only consider every other term in the sequence; see \cite[Section 7]{DHSTmore} for more details.)

For knots in $S^3$, we have the following relationships between $\tau, \varepsilon$, and $\varphi_i$:
\begin{theorem}[{\cite[Proposition 1.2]{DHSTmore}, \cite[Proposition 3.2]{HomEpsilon}}]\label{thm:S3prop}
Let $K$ be a knot in $S^3$. 
\begin{enumerate}
	\item If $\varepsilon(K) = 0$, then $\varphi_i(K) = 0$ for all $i$.
	\item The invariant $\tau$ is equal to
	\[ \tau(K) = \sum_{i=1}^\infty i \varphi_i(K).\]
\end{enumerate}
In particular, if $\varepsilon(K) = 0$, then $\tau(K) = 0$.
\end{theorem}

The invariant $\varepsilon$ can be generalized to a homology concordance invariant \cite[Section 4]{HLL}, which behaves like a sign under connected sum, in the sense that 
\begin{itemize}
	\item if $\varepsilon(Y_1, K_1) =\varepsilon(Y_2, K_2)$, then $ \varepsilon(Y_1 \# Y_2, K_1 \# K_2) = \varepsilon(Y_1,  K_1)$, and 
	\item if $\varepsilon(Y_1, K_1) = 0$, then $\varepsilon(Y_1 \# Y_2, K_1 \# K_2) = \varepsilon(Y_2, K_2)$.
\end{itemize}	
The proof of Theorem \ref{thm:HLL}\eqref{it:HLL2} relies on using the filtered mapping cone of \cite{HeddenLevineSurgery} to produce a manifold-knot pair $(Y, K)$ with $\varepsilon(Y, K) = 0$ and $\tau(Y, K) = 1$. By Theorem \ref{thm:S3prop}, we know that if a knot $K$ in a manifold $Y$ is homology concordant to any knot $J$ in $S^3$ and $\varepsilon (Y,K) = 0$, then $\tau(Y,K) = 0$. Hence $K \subset Y$ is not homology concordant to any knot in $S^3$. Moreover, since $\tau$ is a concordance homomorphism, it follows that any nonzero multiple of $(Y,K)$ has $\varepsilon = 0$ and $\tau$ nonzero, hence cannot be homology concordant to any knot in $S^3$.

\begin{remark}\label{rem:Qslice}
Note that $\varepsilon(K)$, $\tau(K)$, $\varphi_i(K)$, and more generally the local equivalence class of $\CFK(K)$ are all invariant under concordances in rational homology cobordisms. In particular, they all vanish for rationally slice knots. Thus, in order to study $\cC_{\Q S}$, the group of rationally slice knots, we will need additional structure, as discussed in Section \ref{sec:involutiveCFK}.
\end{remark}

\subsection{Ribbon concordances}
As in the case for 3-manifolds, ribbon concordances induce particularly nice maps on the knot Floer complex:
\begin{theorem}[{\cite[Theorem 1.7]{ZemkeRibbon}}]
Let $S$ be a ribbon concordance from $K_0$ to $K_1$ in $S^3 \times [0,1]$, and let $S'$ denote the concordance obtained by reversing $S$. Then 
\[ F_{S'} \circ F_S \co \CFK(K_0) \to \CFK(K_0) \]
is chain homotopic to the identiy, via an $\F[U, V]$-equivariant chain homotopy. In particular, if $S$ is a ribbon concordance from $K_0$ to $K_1$, then $\HFKhat (K_0)$ is a direct summand of $\HFKhat (K_1)$ and $\HFKm(K_0)$ is a direct summand of $\HFKm(K_1)$.
\end{theorem}

\noindent Since knot Floer homology detects the knot genus $g(K)$ \cite{genusbounds}, an immediate consequence of the above theorem is that if there is a ribbon concordance from $K_0$ to $K_1$ then $g(K_0) \leq g(K_1)$ \cite[Theorem 1.3]{ZemkeRibbon}.

\section{Involutive Heegaard Floer homology}\label{sec:HFI}
We would like to use the Heegaard Floer package to study homology cobordism. In light of Theorem \ref{thm:phi}, we see that a richer algebraic structure, namely, chain complexes over a more complicated ring than $\F[U]$, can give us richer invariants. Fortunately for us, Hendricks and Manolescu \cite{HMInvolutive} endowed the Heegaard Floer chain complex with the additional structure of a homotopy involution $\iota$. Very roughly, this additional data lets us think of the Heegaard Floer chain complex as a module over (a quotient of) a two variable polynomial ring, allowing us to employ the techniques used in the proof of Theorem \ref{thm:phi} to define an infinite family of $\Z$-valued homology cobordism homomorphisms. These homomorphisms lead to the proof of Theorem \ref{thm:DHSThomcob}. Furthermore, the characterization of such chain complexes up to a suitable notion of local equivalence is a key ingredient in the proof of Theorem \ref{thm:SFS}.

\subsection{Properties and examples}
Recall that in the construction of Heegaard Floer homology, we specify our 3-manifold $Y$ via a pointed Heegaard diagram $\cH = (\Sigma, \alphas, \betas, z)$, where $\Sigma$ is a closed, oriented surface of genus $g$, and $\alphas$ and $\betas$ are each a collection of $g$ disjoint embedded circles in $\Sigma$. Reversing the orientation of $\Sigma$ reverses the orientation of $Y$, as does reversing the roles of the $\alpha$- and $\beta$-circles. In particular, the Heegaard diagram $\overline{\cH} = (-\Sigma, \betas, \alphas, z)$ describes the same manifold as $\cH$, namely $Y$. Thus, there is a sequence of Heegaard moves taking $\overline{\cH}$ to $\cH$, inducing an $\F[U]$-equivariant chain map
\[ \Phi_{\overline{\cH}, \cH} \co \CFm(\overline{\cH}) \to \CFm(\cH); \] 
this chain map is well-defined since Heegaard Floer homology is natural \cite{JTNaturality}.
There is also a canonical $\F[U]$-equivariant chain complex isomorphism
\[ \eta \co \CFm(\cH) \to \CFm(\overline{\cH}) \]
given by the ``obvious'' identification of the generators of the two chain complexes. Hendricks and Manolescu show that $\iota = \Phi_{\overline{\cH}, \cH} \circ \eta$ is a homotopy involution (that is, that $\iota^2 \simeq  \id$) and prove that for a homology sphere $Y$, the chain homotopy type of the pair $(\CFm(\cH), \iota)$ is an invariant of $Y$ \cite[Proof of Proposotion 2.7]{HMInvolutive}; we will write $(\CFm(Y), \iota)$, called the \emph{$\iota$-complex of $Y$}, to denote a representative of this equivalence class. (An analogous statement holds for a general 3-manifold equipped with a self-conjugate $\spinc$-structure; for ease of exposition, we have chosen to focus on homology spheres to eliminate the need to discuss $\spinc$-structures.) 

\begin{example}
The $\iota$-complex of $S^3$ is $(\F[U], \id)$. (The map $\iota$ is uniquely determined by the fact that $\iota^2 \simeq \id$.)
\end{example}

\begin{example}
The $\iota$-complex of $\Sigma (2, 3, 7)$ is generated over $\F[U]$ by $a, b,$ and $c$ with $\d$, $\iota$, and gradings as follows:
\begin{center}
\begin{tabular}{*{16}{@{\hspace{10pt}}c}}
\hline
&  && $\d$ &&  $\iota$ &&  $\gr$ \\
\hline
& $a$ && $0$ && $c$ && $0$ \\ 
& $b$ && $U a + U c$ &&$b$ &&  $-1$ \\ 
& $c$ && $0$ && $a$ &&  $0$ \\ 
\hline
\end{tabular}
\end{center}
(The map $\iota$ can be computed by realizing $-\Sigma(2,3,7)$ as $+1$-surgery on the left handed trefoil. See \cite[Section 6.8]{HMInvolutive}.)
\end{example}

We have the following K\"unneth-type formula for $\iota$-complexes of connected sums \cite[Theorem 1.1]{HMZConnectedSum}:

\[ (\CFm(Y_1 \# Y_2), \iota_{1\#2}) \simeq (\CFm(Y_1) \otimes_{\F[U]} \CFm(Y_2), \iota_{1} \otimes \iota_{2} ),\]
where $\iota_1$ (respectively $\iota_2$) denotes the homotopy involution on $Y_1$ (respectively $Y_2$) and $\iota_{1 \# 2}$ denotes the homotopy involution on $Y_1 \# Y_2$.

\subsection{Cobordism maps}
The map $\iota$ behaves nicely with respect to cobordism maps. For expositional simplicity, we will focus on homology cobordisms. (One can also consider general cobordisms with a conjugacy class of $\spinc$-structures.) A homology cobordism $W$ from $Y_0$ to $Y_1$ induces a chain map
\[ F_W \co \CFm(Y_1) \to \CFm(Y_2), \]
which commutes with $\iota$, up to homotopy \cite[Proof of Proposition 4.9]{HMInvolutive}:
\[ F_W \circ \iota_1 \simeq \iota_2 \circ F_W. \]
In particular, we can define a refined version of local equivalence as follows:

\begin{definition}\label{def:iotalocequiv}
Two $\iota$-complexes $(C_1, \iota_1)$ and $(C_2, \iota_2)$ are \emph{$\iota$-locally equivalent}, denoted $(C_1, \iota_1) \sim (C_2, \iota_2)$, if there exist $\F[U]$-module chain maps
\[ f \co \CFm(Y_0) \to \CFm(Y_1) \qquad \textup{ and } \qquad g \co \CFm(Y_1) \to \CFm(Y_0) \]
inducing isomorphisms on $U^{-1}\HFm(Y_i)$, such that
\[ f \circ \iota_1 \simeq \iota_2 \circ f \qquad \textup{ and } \qquad g \circ \iota_2 \simeq \iota_1 \circ g. \]
\end{definition}

We can now consider the group
\[  \CFIgroup = \{ (\CFm(Y), \iota) \mid Y \textup{ an integer homology sphere} \} / \sim \]
under the operation induced by tensor product. This construction yields a homomorphism 
\[ \Theta^3_\Z \to  \CFIgroup \]
obtained by sending $[Y]$ to $[(\CFm(Y), \iota)]$. The additional requirement, that the local equivalences homotopy commute with $\iota$, makes this group more interesting than before. However, this group is almost too interesting, in the sense that it is very difficult to understand.

As in the knot case above, a certain algebraic simplification allows us to characterize elements in this group, up to an ever so slightly weaker notion of equivalence. The main idea is to append ``mod $U$'' to every statement in Definition \ref{def:iotalocequiv} involving $\iota$.
\begin{definition}\label{def:almostiotalocequiv}
Two $\iota$-complexes $(C_1, \iota_1)$ and $(C_2, \iota_2)$ are \emph{almost $\iota$-locally equivalent}, if there exist $\F[U]$-module chain maps
\[ f \co \CFm(Y_0) \to \CFm(Y_1) \qquad \textup{ and } \qquad g \co \CFm(Y_1) \to \CFm(Y_0) \]
inducing isomorphisms on $U^{-1}\HFm(Y_i)$, such that
\[ f \circ \iota_1 \simeq \iota_2 \circ f \mod U \qquad \textup{ and } \qquad g \circ \iota_2 \simeq \iota_1 \circ g \mod U. \]
\end{definition}

Similarly, we can relax the definition of an $\iota$-complex so as to only require that $\iota^2 \simeq \id \mod U$; we will call such a complex an \emph{almost $\iota$-complex}.\footnote{Here, we use the word ``almost'' to denote that any statement regarding $\iota$ should be taken mod $U$. In the next section, we use the word ``almost'' to denote that any statement regarding $\iota_K$ should be taken mod $(U, V)$.} Almost $\iota$-local equivalence classes of almost $\iota$-complexes are totally ordered, with
\[ (C_1, \iota_1) \leq (C_2, \iota_2) \]
if there exists an $\F[U]$-module chain map $f \co C_1 \to C_2$ inducing an isomorphism on $H_*(C_i)/U\textup{-torsion}$ 
such that $f \circ \iota_1 \simeq \iota_2 \circ f \mod U$. 

\begin{theorem}[{\cite[Theorem 6.2]{DHSThomcob}}]\label{thm:almostiota}
Every $\iota$-complex is almost $\iota$-locally equivalent to a standard complex, which can be represented by a finite sequence of the form $(a_i, b_i)_{i=1}^n$ where $a_i \in \{\pm 1\}$ and $b_i \in \Z \setminus \{0\}$. Moreover, if we endow the integers with the following unusual order
\[ -1 \lebang -2 \lebang -3 \lebang \dots \lebang 0 \lebang \dots \lebang 3 \lebang 2 \lebang 1 \]
then almost $\iota$-local equivalence classes are ordered lexicographically with respect to their standard representatives.
\end{theorem}

The astute reader may notice that Theorem \ref{thm:almostiota} looks very similar to Theorem \ref{thm:phi}. Indeed, the idea is that $\iota$-complexes can roughly be thought of as chain complexes over $\F[U,Q]/(Q^2)$, which are then very similar to chain complexes over $\F[U, V]$. Analogously, we can define a linearly independent family of homology cobordism homomorphisms
\[ \phi_i \co \Theta^3_\Z \to \Z, \quad i \in \N. \]
These homomorphisms can be used to show that the Brieskorn homology spheres $\Sigma (2j + 1, 4j + 1, 4j + 3)$ span a free infinite rank subgroup of $\Theta^3_\Z$, proving Theorem \ref{thm:DHSThomcob}. Rostovtsev \cite{Rostovtsev} gives an alternate proof of Theorem \ref{thm:almostiota} and extends our result to define an additional, linearly independent integer-valued homology cobordism homomorphism.

Let $\CFIgrouphat$ denote the group of almost $\iota$-complexes up to almost $\iota$-local equivalence. We have the homomorphism
\[ \hat{h} \co  \Theta^3_\Z \to \CFIgrouphat \]
defined by sending $[Y]$ to $[(\CFm(Y), \iota)]$. 
We now describe the main ideas behind the proof of Theorem \ref{thm:SFS}, which states that Seifert fibered spaces do not generate $\Theta^3_\Z$. Let $\Theta_{\SF}$ denote the subgroup of $\Theta^3_\Z$ generated by Seifert fibered spaces. In \cite[Section 8.1]{DHSThomcob}, we determine $\hat{h}(\Theta_\SF)$, the image of $\Theta_\SF$ in $\CFIgrouphat$. Recall that elements of $\CFIgrouphat$ are finite sequences. Elements in $\hat{h}(\Theta_\SF)$ are exactly the sequences satisfying a certain monotonicity condition on their terms; see \cite[Theorem 8.1]{DHSThomcob} for the precise statement. We then use the involutive surgery formula \cite[Theorem 1.6]{HHSZ} to determine the $\iota$-complexes of surgeries on a family of connected sums of torus knots and iterated cables. The sequences associated to these surgeries do not satisfy the monotonicity condition of $\hat{h}(\Theta_\SF)$. A similar statement applies to linear combinations of these surgeries, giving Theorem \ref{thm:SFS}.

\section{Involutive knot Floer homology}\label{sec:involutiveCFK}
In the previous section, we put additional structure, namely the homotopy involution $\iota$, on the Heegaard Floer chain complex $\CFm(Y)$. In this section, we put additional structure, namely a skew-graded, skew-equivariant (i.e., interchanges the actions of $U$ and $V$) chain map $\iota_K$, on the knot Floer chain complex $\CFK(K)$. The map $\iota_K$ is not a homotopy involution; rather, Hendricks-Manolescu \cite[Section 6.2]{HMInvolutive} show that $\iota_K$ squares to be homotopic to the Sarkar map \cite{SarkarMovingBasepoints}, which is induced by moving the basepoints $w$ and $z$ once around the knot $K$.
The map $\iota_K$ is the additional structure alluded to in Remark \ref{rem:Qslice}.

\subsection{Properties and examples}
The map $\iota_K$ is defined in a similar way to $\iota$. Consider a doubly pointed Heegaard diagram $\cH = (\Sigma, \alphas, \betas, z, w)$ for a knot $K$ in an integer homology sphere $Y$. (With minor modifications, these constructions also work for null-homologous knots in any 3-manifold.) The doubly pointed Heegaard diagram $\overline{\cH} = (-\Sigma, \betas, \alphas, w, z)$ also describes $K \subset Y$, and thus there is a sequence of Heegaard moves from $\overline{\cH}$ to $\cH$ inducing a $\F[U, V]$-equivariant chain map
\[ \Phi_{\overline{\cH}, \cH} \co \CFK(\overline{\cH}) \to \CFK(\cH). \]
(Care should be taken with regards to the basepoints; see \cite[Section 6.2]{HMInvolutive} for details.)
There is also a canonical skew-equivariant isomorphism
\[ \eta_K \co \CFK(\cH) \to \CFK(\overline{\cH}), \]
given by the ``obvious'' identification of the generators of the two chain complexes. Then $\iota_K$ is defined to be $\Phi_{\overline{\cH}, \cH} \circ \eta_K$. The chain homotopy type of the pair $(\CFK(\cH), \iota_K)$ is an invariant of the knot $K$ in $Y$ \cite[Proposition 6.3]{HMInvolutive}. As usual, we write $(\CFK(K), \iota_K)$ to denote a representative of this chain homotopy equivalence class, and we call the pair $(\CFK(K), \iota_K)$ an \emph{$\iota_K$-complex}.

\begin{example}
The $\iota_K$-complex of the unknot in $S^3$ is $(\F[U, V], \id)$. (The map $\iota_K$ is uniquely determined by the fact that it squares to be homotopic to the Sarkar map.)
\end{example}

\begin{example}
The $\iota_K$-complex of the right handed trefoil is described by:
\begin{center}
\begin{tabular}{*{16}{@{\hspace{10pt}}c}}
\hline
&  && $\d$ && $\iota_K$ && $\gr$ \\
\hline
& $a$ && $0$ && $c$ &&  $(0, -2$) \\ 
& $b$ && $U a + V c$ && $b$ &&  $(-1, -1)$ \\ 
& $c$ && $0$ && $a$ &&  $(-2, 0)$ \\ 
\hline
\end{tabular}
\end{center}
(The map $\iota_K$ is uniquely determined by the fact that it is skew-graded and squares to be homotopic to the Sarkar map.)
\end{example}

There is a K\"unneth-type formula for $\iota_K$-complexes \cite[Theorem 1.1]{ZemConnectedSums}:
\[ (\CFK( Y_1 \# Y_2, K_1 \# K_2), \iota_{K_1 \# K_2} )\simeq (\CFK(K_1) \otimes \CFK(K_2),  \iota_{K_1} \otimes \iota_{K_2} + (\Phi \otimes \Psi ) \circ (\iota_{K_1} \otimes \iota_{K_2})). \]
See \cite[Section 4.2]{ZemConnectedSums} for the definitions of $\Phi$ and $\Psi$; for an expository overview, see \cite[Section 2]{HKPS}.

\subsection{Maps induced by concordances}
Let $(W,S)$ be a $\Z/2\Z$-homology cobordism $(Y_0, K_0)$ to $(Y_1, K_1)$; that is, $W$ is a $\Z/2\Z$-homology cobordism from $Y_0$ to $Y_1$ and $S$ is a concordance from $K_0$ to $K_1$. (More generally, one can consider spin cobordisms; see \cite{HMInvolutive}.)  As one would hope, the module homomorphism $F_{W,S}$ induced by $(W, S)$ behaves nicely with respect to $\iota_K$ \cite[Theorem 1.3]{ZemConnectedSums} in the sense that 
\[ F_{W,S} \circ \iota_{K_0} \simeq \iota_{K_1} \circ F_{W, S}. \]

Based on previous sections, it may now be apparent to the reader what we do next. We jump straight to the definition of almost $\iota_K$-local equivalence; the definition of $\iota_K$-local equivalence can be obtained by striking out both instances of ``mod $(U, V)$'' in the definition below.

\begin{definition}
Two $\iota_K$-complexes $(C_1, \iota_1)$ and $(C_2, \iota_2)$ are \emph{almost $\iota_K$-locally equivalent} if there exist $\F[U, V]$-module chain maps
\[ f \co \CFK(Y_0, K_0) \to \CFK(Y_1, K_1) \quad \textup{ and } \quad g \co \CFK(Y_1, K_1) \to \CFK(Y_0, K_0) \]
inducing isomorphisms on $(U, V)^{-1}\HFK(Y_i, K_i)$ such that
\[ f \circ \iota_{K_0} \skewsimeq \iota_{K_1} \circ f \mod (U, V) \quad \textup{ and } \quad  g \circ \iota_{K_1} \skewsimeq \iota_{K_0} \circ g \mod (U, V),  \]
where $\skewsimeq$ denotes skew-equivariant homotopy equivalence.
\end{definition}

We can now consider $\CFKIgrouphat$, the group of $\iota_K$-complexes modulo almost $\iota_K$-local equivalence, with the operation induced by tensor product. Note that this group has 2-torsion, generated by, for example, the figure-eight knot, and hence this group is not totally ordered. In particular, there are rationally slice knots, such as the figure-eight, with nontrivial image in $\CFKIgrouphat$; this is possible because the (almost) local equivalence class of a knot is an invariant of concordances in $\Z/2\Z$-homology cobordisms, rather than $\Q$-homology cobordisms, as is the case in the non-involutive setting. Thus, this toolkit is particularly well-equipped for studying rationally slice knots.

The proof of Theorem \ref{thm:Qslice} relies on finding a linearly independent family of rationally slice knots. The knots under  consideration are $K_n$, the $(2n+1, 1)$-cable of the figure-eight; these knots are rationally slice because the figure-eight is rationally slice, and thus its $(2n+1, 1)$-cable is rationally concordant (i.e., concordant in a rational homology cobordism) to the $(2n+1, 1)$-torus knot, which is the unknot.

We compute the almost $\iota_K$-local equivalence class of $(\CFK(K_n), \iota_K)$ using bordered Floer homology \cite{LOTBordered, HRW}, in particular, its applications to cables \cite{HW-cables, Petkova}, together with formal properties of $\iota_K$, such as the fact that it squares to be homotopic to the Sarkar map. With this computation in hand, we observe that
 there is certain structure in $(\CFK(K_n), \iota_K)$ (roughly, a particular $\F[U]/U^n$ summand in $\HFKm(K_n)$) which, using the formula for connected sums of $\iota_K$-complexes and properties of almost $\iota_K$-local equivalences, allows us to determine that the $K_n$ are linearly independent.

\section{What next?}\label{sec:future}
As we've demonstrated, the Heegaard Floer package can answer a range of questions in low-dimensional topology. Do these techniques have the potential to answer any of the open questions from Section \ref{sec:intro}?

Question \ref{quest:thetatorsion} asks whether $\Theta^3_\Z$ contains any torsion. The most likely torsion is order two, generated by a homology sphere admitting an orientation-reversing self-diffeomorphism. There are many constructions for building homology spheres with orientation-reversing self-diffeomorphisms (for example, the double-branched cover of an amphichiral knot with determinant one, or the splice of a knot complement with that of its mirror), but so far, there has been no success in obstructing  such an example from being homology cobordant to $S^3$. One can consider an algebraic version of the local equivalence group $\CFIgroup$, by considering all $\iota$-complexes (not just those known to be realized by a 3-manifold) modulo local equivalence. This algebraic group is known to have two-torsion; the difficulty lies in finding a 3-manifold that realizes such an algebraic example. At present, computations of $\iota$-complexes are limited to certain special families of manifolds (e.g., Seifert fibered spaces, surgeries on knots in $S^3$); we hope to improve this shortcoming in the future.

Question \ref{quest:surgeries} asks whether surgeries on knots in $S^3$ generate $\Theta^3_\Z$. This seems like a hard question to answer with Heegaard Floer homology, as the question about which $\iota$-complexes can be realized by surgery on a knot in $S^3$ then reduces to the question of which $\iota_K$-complexes can be realized by knots in $S^3$. Even without the additional structure of $\iota_K$, this is a difficult question; for some partial answers, see \cite{HeddenWatson}, as well as  more recent progress in \cite{BaldwinVV}, \cite{Krcatovichred}, and \cite{NiNtT}.
 
As for Questions \ref{quest:ZHBslice} and \ref{quest:S4PC}, which ask for knots that are not slice in $B^4$ but are slice in a homology $B^4$ or a homotopy $B^4$, respectively, it seems unlikely that Heegaard Floer homology will be able to provide an answer. Indeed, with the current technology, if the Heegaard Floer package obstructs a knot from being slice in $B^4$, then it also obstructs the knot from being slice in a homology or homotopy $B^4$. There are other invariants which may be able to shed light on this question. For example, at present, it remains open whether or not the Rasmussen $s$-invariant \cite{RasmussenS}, defined using the Lee \cite{LeeKh} deformation of Khovanov homology \cite{Khovanov} (see \cite{BarNatan} for an expository overview), vanishes for knots that are slice in a homology or homotopy $B^4$.

We now turn to Question \ref{quest:QHBslice}, which asks whether there exists a knot $K \subset S^3$ that is not slice in $B^4$ but is slice in a rational homology 4-ball $W$ with $|H_1(W; \Z)|$ odd. It seems unlikely that Heegaard Floer homology, in its present form, can address this condition. Note that involutive Heegaard Floer homology gives obstructions to being slice in a $\Z/2\Z$-homology 4-ball; if $W$ is a rational homology 4-ball with $|H_1(W; \Z)|$ odd, then it is a $\Z/2\Z$-homology ball. However, recall that prior to the advent of involutive Heegaard Floer homology, there was no way to use Heegaard Floer homology to obstruct a knot (such as the figure-eight) from being slice in any rational homology 4-ball. Perhaps there is some other additional structure that we can add to the Heegaard Floer package, yielding new obstructions. Alternatively, it remains possible that the $s$-invariant may have something to say about this question. 

Conjecture \ref{conj:GordonRibbon} posits that ribbon concordance is a partial order. Zemke \cite[Theorem 1.7]{ZemkeRibbon} proved that if there is a ribbon concordance from $K_0$ to $K_1$, then $\HFKhat(K_0)$ injects into $\HFKhat(K_1)$. Thus, if there is also a ribbon concordance from $K_1$ to $K_0$, then $\HFKhat(K_0) \cong \HFKhat(K_1)$. Note that there are infinite families of knots with the same knot Floer homology \cite[Theorem 1]{HeddenWatson}. However, as far as the author knows, there are no known ribbon concordances between distinct knots in any of those families. Further investigation is needed before we rule out knot Floer homology as a tool for resolving Conjecture \ref{conj:GordonRibbon}.


Closely related is Conjecture \ref{conj:ribbon3}, which posits that ribbon $\Q$-homology cobordism is a partial order on 3-manifolds. There is a ribbon homology cobordism from $S^3$ to $Y \# -Y$ for any homology sphere $Y$. Taking $Y= \Sigma(2,3, 5)$, and noting that $\HFm(\Sigma(2,3,5) \# -\Sigma(2,3,5)) \cong \HFm(S^3) \cong \F_{(0)}[U]$, we see that we have two distinct 3-manifolds with the same Heegaard Floer homology and a ribbon homology cobordism in one direction. (As alluded to above, we do not know of an analogous example in the knot case.) However, since $\Sigma(2,3,5) \# -\Sigma(2,3,5)$ does not bound a simply-connected homology 4-ball \cite[Proposition 1.7]{TaubesPeriodic}, it follows that there is no ribbon homology cobordism from $\Sigma(2,3,5) \# -\Sigma(2,3,5)$ to $S^3$ (for if there were, we could glue a 4-ball to the $S^3$ end and obtain a simply-connected homology ball with boundary $\Sigma(2,3,5) \# -\Sigma(2,3,5)$ ). We refer the reader to \cite{DLVW} for further evidence, some of it coming from various Floer homologies, in support of Conjecture \ref{conj:ribbon3}.

As we have seen, advances in Heegaard Floer homology have answered many questions about homology cobordism and knot concordance. These successes were not immediate; they began in 2003, when Ozsv\'ath-Szab\'o \cite{OSIntersectionForms, OS4ballgenus}  defined the homomorphisms
\[ d \co \Theta^3_\Z \to 2\Z \qquad \textup{ and } \qquad \tau \co \cC \to \Z. \]
The next major step in extracting concordance information from the knot Floer complex was the definition of $\varepsilon$ in 2014 \cite{HomEpsilon}, which in turn led to two infinite families of concordance homomorphisms: 
\[ \Upsilon_t \co \cC \to \R, \ t \in [0, 2] \qquad \textup{ and } \qquad \varphi_i \co \cC \to \Z, \ i \in \N, \]
the former defined by Ozsv\'ath-Stipsicz-Szab\'o's \cite{OSSUpsilon}, and the latter by Dai, Stoffregen, Truong, and the author \cite{DHSTmore}. The algebraic framework necessary to define $\Upsilon_t$ and $\varphi_i$ existed since the inception of knot Floer homology in the early 2000s, yet it took over a decade for anyone to exploit this structure to define these new homomorphisms. Concurrent with these developments was the advent of Hendricks-Manolescu's involutive Heegaard Floer homology \cite{HMInvolutive}, which put new, more refined structure on the Heegaard Floer and knot Floer homology packages, yielding new homology cobordism homomorphisms and new rational concordance obstructions. We look forward to seeing whether the Heegaard Floer package in its present form can be further mined for new applications, to refining the structure on these invariants even more to prove new theorems, and to developing new, unanticipated tools for resolving the questions and conjectures that we have posed here.


\bibliographystyle{alpha}
\bibliography{biblio}

\end{document}